\documentclass[12pt]{amsart}
\usepackage{amssymb}
\newtheorem{theorem}{Theorem}[section]
\newtheorem{lemma}[theorem]{Lemma}
\theoremstyle{definition}
\newtheorem{definition}{Definition}[section]

\newtheorem{example}{Example}[section]

\theoremstyle{remark}

\numberwithin{equation}{section}

\begin{document}
\title[On generalized $\phi$-recurrent generalized $(k, \mu)$-contact metric manifolds]{On generalized $\phi$-recurrent generalized $(k, \mu)$-contact metric manifolds}
\author[S. K. Hui]{Shyamal Kumar Hui}
\begin{abstract}
The present paper deals with the study of generalized $\phi$-recurrent generalized $(k,\mu)$-contact metric 
manifolds with the existence of such notion by a proper example.
\end{abstract}
\maketitle
\section{Introduction}
In 1995 Blair, Koufogiorgos and Papantoniou \cite{BLAIR3} introduced the notion of $(k,\mu)$-contact metric manifolds, where $k$ and $\mu$ are real constants and a full classification of such manifolds was given by Boeckx \cite{BOEC}. Assuming $k, \mu$ be smooth functions, Koufogiorgos and Tsichlias \cite{KOUF} introduced the notion of generalized
$(k, \mu)$-contact metric manifolds with the existence of such notions.\\
\indent The notion of local symmetry of a Riemannian manifold has been weakened by many authors in several ways to a different extent. As a weaker version of local symmetry, Takahashi \cite{TAKA} introduced the notion of local $\phi$-symmetry on a Sasakian manifold. Recently Shaikh \cite{SHAIKH2} studied the locally $\phi$-symmetry generalized $(k,\mu)$-contact metric manifolds. Also Baishya, Eyasmin and Shaikh \cite{BAISHYA} introduced and studied the locally $\phi$-recurrent $(k,\mu)$-contact metric manifolds and locally $\phi$-recurrent generalized $(k,\mu)$-contact metric manifolds. Generalizing all these notions of local $\phi$-symmetry, in the present paper we introduce generalized $\phi$-recurrent
generalized $(k,\mu)$-contact metric manifolds.\\
\indent In 1979, Dubey \cite{DUBEY} introduced generalized recurrent manifolds. We note that generalized recurrent manifolds are also studied in
(\cite{ARSL}, \cite{DE}). A Riemannian manifold (M,g) is called generalized recurrent \cite{DE} if its curvature tensor $R$ satisfies the condition
\begin{equation}\label{eqn1.1}
  \nabla R=A \otimes R+B \otimes G
\end{equation}
where $A$ and $B$ are two non-vanishing $1$-forms defined by $A(.)=g(.~,\rho _1)$, $B(.)=g(.~,\rho _2)$ and the tensor $G$ is defined by
\begin{equation}\label{eqn1.2}
  G(X,Y)Z=g(Y,Z)X-g(X,Z)Y
\end{equation}
for all $X,Y,Z \in \chi(M)$; $\chi(M)$ being the Lie algebra of the smooth vector fields and $\nabla$ denotes covariant differentiation
with respect to the metric $g$. Here $\rho _1$ and $\rho_2$ are vector fields associated with $1$-forms $A$ and $B$ respectively.\\
Especially, if the $1$-form $B$ vanishes, then (\ref{eqn1.1}) turns into the notion of recurrent manifold introduced by Walker \cite{WALK}.\\
 \indent A Riemannian manifold $(M,g)$ is called a generalized Ricci-recurrent \cite{DE1} if its Ricci tensor $S$ of type $(0,2)$ satisfies
 the condition
 \begin{equation}\label{eqn1.3}
   \nabla S=A \otimes S+B\otimes g
 \end{equation}
where $A$ and $B$ are defined in (\ref{eqn1.1}).\\
\indent In particular, if $B=0$, then (\ref{eqn1.3}) reduces to the notion of Ricci-recurrent manifolds introduced by Patterson \cite{PATT}.\\
\indent Recently Shaikh and Ahmad \cite{SHAIKH3}
 introduced the notion of generalized $\phi$-recurrent Sasakian manifolds. The present paper deals with the study of generalized $\phi$-recurrent generalized $(k,\mu)$-contact metric manifolds. The paper is organized as follows. Section $2$ is concerned with some preliminaries.
In section $3$, we study generalized $\phi$-recurrent generalized $(k,\mu)$-contact metric manifolds. Finally, we construct an example of a generalized $\phi$-recurrent generalized $(k,\mu)$-contact metric manifold which is neither $\phi$-symmetric nor $\phi$-recurrent in the last section.
\section{Prelminaries}
A contact manifold is a $C^\infty$ manifold $M^{2n+1}$ equipped with a global $1$-form $\eta$ such that $\eta \wedge(d \eta)^n \neq0$
everywhere on $M^{2n+1}$. Given a contact form $\eta$ it is well known that there exists a unique vector field $\xi$, called the characteristic
vector field of $\eta$, such that $\eta(\xi)=1$ and $d \eta(X,\xi)=0$ for every vector field $X$ on $M^{2n+1}$. A Riemannian metric is said to be associated
metric if there exists a tensor field $\phi$ of type $(1,1)$ such that
\begin{equation}\label{eqn2.1}
   d \eta(X,Y)=g(X,\phi Y), \ \ \ \ \ \ \ \eta(X)=g(X,\xi),
\end{equation}
\begin{equation}\label{eqn2.2}
  \phi \xi=0,  \ \ \ \ \ \eta(\phi X)=0,   \ \ \ \ \ \ \phi^2 X=-X+\eta(X)\xi,
\end{equation}
\begin{equation}\label{eqn2.3}
  g(\phi X, \phi Y)=g(X,Y)-\eta(X) \eta(Y)
\end{equation}
for all vector fields $X,Y$ on $M^ {2n+1}$. Then the structure $(\phi, \xi, \eta, g)$ on $M^{2n+1}$ is called a contact metric structure and the
manifold $M^{2n+1}$ equipped with such structure is called a contact metric manifold \cite{BLAIR1}.\\\
\indent Given a contact metric manifold $M^{2n+1}(\phi, \xi, \eta, g)$ we define a $(1,1)$ tensor field $h$ by $h=\frac{1}{2} \pounds_\xi \phi,$ where
$\pounds$ denotes the Lie differentiation. Then $h$ is symmetric and satisfies $h \phi=-\phi h$. Thus, if $\lambda$ is an eigenvalue of $h$ with eigenvector $X$, $-\lambda$ is also an eigenvalue with eigenvector $\phi X$. Also we have $Tr.~h=Tr.~\phi h=0$ and $h \xi=0.$ Moreover, if $\nabla$ denotes the Riemannian connection of $g$, then the following relation holds:
\begin{equation}\label{eqn2.4}
  \nabla _X \xi=-\phi X-\phi hX.
\end{equation}
The vector field $\xi$ is Killing vector with respect to $g$ if and only if $h=0.$
A contact metric manifold $M^{2n+1}(\phi, \xi, \eta,g)$ for which $\xi$ is a Killing vector is said to be a $K$-contact manifold. A contact structure on $M^{2n+1}$ gives rise to an almost complex structure on the product $M^{2n+1}\times R$. If this almost complex structure is integrable, the contact metric manifold is said to be Sasakian. Equivalently, a contact metric manifold
is Sasakian if and only if the relation
\begin{equation*}
  R(X,Y)\xi=\eta(Y)X-\eta(X)Y
\end{equation*}
holds for all $X,Y$, where $R$ denotes the curvature tensor of the manifold.
\begin{lemma}\cite{BLAIR2}
   Let $M^{2n+1}(\phi,\xi,\eta,g)$ be a contact metric manifold with $R(X,Y)\xi=0$ for all vector fields $X,Y$ tangent to $M$. Then
   M is locally isometric to the Riemannian product $E^{n+1}(0)\times S^n(4)$.
 \end{lemma}
For a contact metric manifold $M^{2n+1}(\phi, \xi, \eta, g)$, the $(k, \mu)$-nullity distribution is
\begin{eqnarray*}
p\rightarrow N_{p}(k,\mu)&=& [Z \in {T_p} M: R(X,Y)Z=k\{(g(Y,Z)X-g(X,Z)Y\}\\
\nonumber&+& \mu \{g(Y,Z)hX-g(X,Z)hY\}]
\end{eqnarray*}
 for any $X,Y \in T_pM$, $k, \mu$ are real numbers. Hence, if the characteristic vector field $\xi$ belongs to the $(k, \mu)$-nullity
 distribution, then we have
 \begin{equation}\label{eqn3.1}
   R(X,Y)\xi=k[\eta(Y)X-\eta(X)Y]+\mu [\eta(Y)hX-\eta(X)hY].
 \end{equation}
 Thus a contact metric manifold satisfying relation (\ref{eqn3.1}) is called a $(k, \mu)$-contact metric manifold \cite{BLAIR3}.
 In particular, if $\mu=0$, then the notion of $(k, \mu)$-nullity distribution reduces to the notion of $k$-nullity distribution
 reduces to the notion of $k$-nullity distribution, introduced by Tanno \cite{TANNO}. A $(k,\mu)$-contact metric manifold is Sasakian
 if and only if $k=1$. In a $(k, \mu)$-contact metric manifold the following relations hold (\cite{BLAIR3}, \cite{SHAIKH1}):
 \begin{equation}\label{eqn3.2}
   h^2=(k-1)\phi^2,\ \ \ \ \ k\leq 1,
 \end{equation}
  \begin{equation}\label{eqn3.3}
    (\nabla_X \phi)(Y)=g(X+hX,Y)\xi-\eta(Y)(X+hX),
  \end{equation}
  \begin{eqnarray}
  (\nabla_ X h)(Y) &=& \{(1-k)g(X,\phi Y)+g(X,h\phi Y)\}\xi \\
  \nonumber &+& \eta(Y)[h(\phi X+\phi hX)]-\mu{\eta(X)}\phi h Y,
  \end{eqnarray}
  \begin{equation}\label{eqn3.5}
    R(\xi,X)Y=k[g(X,Y)\xi-\eta(Y)X]+\mu [g(hX,Y)\xi-\eta(Y)hX],
  \end{equation}
 \begin{eqnarray}
  \label{eqn3.6}
    \eta(R(X,Y)Z) &=& k[g(Y,Z)\eta(X)-g(X,Z)\eta(Y)]\\
    \nonumber &+& \mu[g(h Y,Z)\eta(X)-g(h X, Z)\eta(Y)],
 \end{eqnarray}
 \begin{equation}\label{eqn3.7}
   S(X, \xi)=2nk \eta(X),
 \end{equation}
 \begin{equation}\label{eqn3.8}
   Q \phi-\phi Q=2[2(n-1)+\mu]h \phi,
 \end{equation}
 \begin{eqnarray}
 \label{eqn3.9}
 S(X,Y) &=& [2(n-1)-n \mu]g(X,Y)+[2(n-1)+\mu]g(hX,Y)\\
\nonumber &+& [2(1-n)+n(2k+\mu)]{\eta(X)\eta(Y)},  \ \ \ \ \ {n\geq 1,}
 \end{eqnarray}
 \begin{equation}\label{eqn3.10}
   r=2n(2n-2+k-n \mu),
 \end{equation}
 \begin{equation}\label{eqn3.11}
   S(\phi X,\phi Y)=S(X,Y)-2nk \eta(X) \eta(Y)-2(2n-2+\mu)g(hX,Y),
 \end{equation}
 where $S$ is the Ricci tensor of type $(0,2)$, $Q$ is the Ricci-operator, i.e., $g(QX,Y)=S(X,Y)$ and $r$ is the scalar curvature of the manifold.
 From (\ref{eqn2.4}), it follows that
 \begin{equation}\label{eqn3.12}
   (\nabla_X \eta)(Y)=g(X+hX,\phi).
 \end{equation}
 Also we have from (\ref{eqn3.1}) that
 \begin{eqnarray}
 \label{eqn3.13}
   (\nabla_W R)(X,Y)\xi &=& k[g(W+hW,\phi Y)X-g(W+hW,\phi X)Y]\\
   \nonumber &+& \mu[g(W+hW,\phi Y)hX-g(W+hW,\phi X)hY \\
   \nonumber&+&\{(1-k)g(W,\phi X)+g(W,h\phi X)\}\eta(Y)\xi\\
   \nonumber&-& \{(1-k)g(W,\phi Y)+g(W,h\phi Y)\}\eta(X)\xi\\
   \nonumber&+& \mu \eta(W)\{\eta(X)\phi hY-\eta(Y)\phi hX\}]\\
   \nonumber&+& R(X,Y)\phi W+R(X,Y)\phi hW.
 \end{eqnarray}

 \section{generalized $\phi$-recurrent $(k,\mu)$-contact metric manifolds}
 \begin{definition}
   A generalized $(k, \mu)$-contact metric manifold $(M^n,g)$ is said to be a generalized $\phi$-recurrent generalized $(k,\mu)$-contact
   metric manifold if the relation
   \begin{equation}\label{eqn4.1}
     \phi^2((\nabla_W R)(X,Y)Z)=A(W)\phi^2(R(X,Y)Z)+B(W)\phi^2(G(X,Y)Z)
   \end{equation}
   holds for all $X,Y,Z,W \in \chi(M)$ and $A$ and $B$ are two non vanishing $1$-forms such that $A(X)=g(X, \rho_1)$,
   $B(X)=g(X,\rho_2).$ Here $\rho_1$ and $\rho_2$ are vector fields associated with $1$-forms $A$ and $B$ respectively.
 \end{definition}
 Let us consider a generalized $\phi$-recurrent generalized $(k,\mu)$-contact metric manifold. Then by virtue of (\ref{eqn2.2}),
 we have from (\ref{eqn4.1}) that
\begin{align}\label{eqn4.2}
 &-(\nabla_W R)(X,Y)Z+\eta((\nabla_W R)(X,Y)Z)\xi &\\
\nonumber& = A(W)[-R(X,Y)Z+\eta(R(X,Y)Z)\xi]&\\
\nonumber&+ B(W)[-G(X,Y)Z+\eta(G(X,Y)Z)\xi]
\end{align}
from which it follows that
 \begin{flushleft}
\begin{align}\label{eqn4.3}
-g((\nabla_W R)(X,Y)Z,U)+\eta((\nabla_W R)(X,Y)Z)\eta(U) &  \\
\nonumber  = A(W)[-g(R(X,Y)Z,U)+\eta(R(X,Y)Z)\eta(U)] & \\
\nonumber +B(W)[-g(G(X,Y)Z,U)+\eta(G(X,Y)Z)\eta(U)].
\end{align}
\end{flushleft}
Taking an orthonormal frame field and then contracting (\ref{eqn4.3}) over $X$ and $U$ and then using (\ref{eqn1.2}) and (\ref{eqn3.5}), we get
\begin{align}
\label{eqn4.4}
&-(\nabla_W S)(Y,Z)+g((\nabla_W R)(\xi, Y)Z,\xi)&\\
\nonumber&= A(W)[-S(Y,Z)+k\{g(Y,Z)-\eta(Y) \eta(Z)\}+\mu\{g(hY,Z)-\eta(Z)\eta(hY)\}]\\
\nonumber&+B(W)[-(2n-1)g(Y,Z)-\eta(Y)\eta(Z)].
\end{align}
Plugging $Z=\xi  $ in (\ref{eqn4.4}), we obtain
\begin{equation}\label{eqn4.5}
  (\nabla_W S)(Y,\xi)=A(W)S(Y,\xi)+2nB(W)\eta(Y).
\end{equation}
By virtue of (\ref{eqn2.4}), (\ref{eqn3.7}) and (\ref{eqn3.12}) it follows from (\ref{eqn4.3}) that
\begin{align}\label{eqn4.6}
  2nkg(W+hW,\phi Y)+S(Y,\phi W+\phi h W) &=2n[kA(W)+B(W)]\eta(Y).
\end{align}
Setting $Y=\xi$ in (\ref{eqn4.6}) and using (\ref{eqn2.2}) and (\ref{eqn3.7}), we get
\begin{equation}\label{eqn4.7}
  k A(W)+B(W)=0.
\end{equation}
In view of (\ref{eqn4.7}), (\ref{eqn4.6}) yields
\begin{equation}\label{eqn4.8}
  S(Y,\phi W+\phi h W)=2nkg(Y,\phi W+\phi hW).
\end{equation}
Replacing $Y$ by $\phi Y$ in (\ref{eqn4.8}) and using (\ref{eqn2.3}) and (\ref{eqn3.11}), we get
\begin{equation}\label{eqn4.9}
  S(Y, W+hW)=2nkg(Y,W+hW)+2(2n-2+\mu)g(hY,W+hW).
\end{equation}
Again replacing $Y$ by $hY$ in (\ref{eqn4.9}) and  using (\ref{eqn2.2}) and (\ref{eqn3.2}), we get

\begin{eqnarray}
\label{eqn4.9a}
 \ \ \  S(Y,hW)-(k-1)S(Y,W)&=&-2nk(k-1)g(Y,W)\\
 \nonumber &-&2(k-1)(2n-2+\mu)g(h Y,W)\\
  \nonumber &+&2(k-1)(2n-2+\mu)\eta(W)\eta(hY).
\end{eqnarray}
Subtracting (\ref{eqn4.9a}) from (\ref{eqn4.9}), we get
\begin{eqnarray}
\label{eqn4.9b}
 kS(Y,Z)&=&2nk^2g(Y,W)+2k(2n-2+\mu)g(hY,W)\\
 \nonumber&-&2(k-1)(2n-2+\mu)\eta(W)\eta(hY).
\end{eqnarray}
This leads to the following:
\begin{theorem}
  In a generalized $\phi$-recurrent generalized $(k,\mu)$-contact metric manifold, the $1$-forms $A$ and $B$ are related by the relation $(\ref{eqn4.7})$
  and the Ricci tensor $S$ is of the form $(\ref{eqn4.9b})$.
  \end{theorem}
\indent Changing $W,X,Y$ cyclically in (\ref{eqn4.3}) and adding them we get by virtue of Bianchi identity and using (\ref{eqn4.7}), we get
\begin{align}\label{eqn4.10}
&A(W)[-g(R(X,Y)Z,U)+kg(G(X,Y)Z,U)+\{\eta(R(X,Y)Z)&\\
\nonumber&-k \eta(G(X,Y)Z)\}\eta(U)]+A(X)[-g((R(Y,W)Z,U)+kg(G(Y,W)Z,U)&\\
\nonumber&+\{\eta(R(Y,W)Z)-k\eta(G(Y,W)Z)\eta(U)\}]+A(Y)[-g(R(W,X)Z,U)&\\
\nonumber&+kg(G(W,X)Z),U)+\{\eta(R(W,X)Z)- k{\eta(G(W,X)Z)}\}\eta(U)]=0.
\end{align}
Contracting (\ref{eqn4.10}) over $Y$ and $Z$, we get
\begin{align}\label{eqn4.11}
  &A(W)[-S(X,U)+2nkg(X,U)]-A(X)[-S(W,U)+2nkg(W,U)]&\\
\nonumber&+A(R(W,X)U)+k\{A(X)g(W,U)-A(W)g(X,U)\}-A(R(W,X)\xi)\eta(U)&\\
\nonumber&-k\{A(X)\eta(W)-A(W)\eta(X)\}\eta(U)=0.
\end{align}
Again contracting (\ref{eqn4.11}) over $X$ and $U$ and using (\ref{eqn3.1}), we get
\begin{equation}\label{eqn4.12}
  2A(QW)-[r-2n(2n-1)]A(W)-\mu A(hW)=0.
\end{equation}
This leads to the following:
\begin{theorem}
  In a generalized $\phi$-recurrent generalized $(k,\mu)$-contact metric manifold, the relation $(\ref{eqn4.12})$ holds for all $W$.
\end{theorem}
Using (\ref{eqn3.5}), (\ref{eqn3.13}) and the relation $g((\nabla_W R)(X,Y)Z,U)=-g((\nabla_W R)(X,Y)U,Z)$, we have
\begin{eqnarray}
\label{eqn4.13}
  g((\nabla_W R)(\xi,Y)Z,\xi) &=& \mu[\{(1-k)g(W, \phi Y)+g(W,h \phi Y)\\
  \nonumber&-&g(hY,\phi(W+hW))\}\eta(Z)-\mu \eta(W)g(\phi hY,Z)].
\end{eqnarray}
By virtue of (\ref{eqn4.7}) and (\ref{eqn4.13}) it follows from (\ref{eqn4.4}) that
\begin{eqnarray}
 \label{eqn4.14}
  (\nabla_W S)(Y,Z) &=& A(W)S(Y,Z)-2nkA(W)g(Y,Z)\\
  \nonumber &+& \mu[\{A(W)\eta(h Y)-(1-k)g(W,\phi Y)\\
  \nonumber&-& g(W,h\phi Y)+g(h Y,\phi(W+h W))\}\eta(Z)\\
  \nonumber&-& A(W)g(h Y,Z)+\mu\eta(W)g(\phi h Y,Z)].
\end{eqnarray}
This leads the following:
\begin{theorem}
  A generalized $\phi$-recurrent generalized $(k,\mu)$-contact metric manifold is generalized Ricci recurrent if and only if the following relation holds:
\begin{align*}
&\{A(W)\eta(h Y)-(1-k)g(W,\phi Y)-g(W,h \phi Y)+g(hY,\phi(W+hW))\}\eta(Z)&\\
\nonumber&-A(W)g(hY, Z)+\mu\eta(W)g(\phi hY,Z)=0.
\end{align*}
\end{theorem}
Again from (\ref{eqn4.2}), we get
\begin{eqnarray}
\label{eqn4.15}
\ \ (\nabla_W R)(X,Y)\xi &=& A(W)[k\{\eta(Y)X-\eta(X)Y\}+\mu \{\eta(Y)hX-\eta(X)hY\}]\\
  \nonumber&+& B(W)[\eta(Y)X-\eta(X)Y].
\end{eqnarray}
From (\ref{eqn3.13}) and (\ref{eqn4.15}) we obtain
\begin{align}\label{eqn4.16}
&k{\big[}g(W+hW,\phi Y)X-g(W+hW, \phi X)Y{\big]}&  \\
  \nonumber &+\mu {\big[}g(W+hW, \phi Y)hX-g(W+hW, \phi X)hY+\{(1-k)g(W,\phi X)\\
 \nonumber &+g(W,h \phi X)\eta(Y)\xi -\{(1-k)g(W,\phi Y)+g(W, h \phi Y)\}\eta(X)\xi &  \\
 \nonumber &+\mu \eta(W)\{\eta(X)\phi hY-\eta(Y)\phi hX\}]+R(X,Y)\phi W &  \\
  \nonumber&+R(X,Y)\phi h W-B(W)[\eta(Y)X-\eta(X)Y{\big]} &  \\
  \nonumber&-A(W){\big[}k\{\eta(Y)X-\eta(X)Y\}+\mu \{\eta(Y)h X-\eta(X)h Y\}{\big]}=0.
\end{align}
Replacing $W$ by $\phi W$  in (\ref{eqn4.16}) and using (\ref{eqn2.2}) we get
\begin{align}\label{eqn4.17}
&R(X,Y)W+R(X,Y)hW &\\
\nonumber&=k{\big[}g(W+hW,Y)hX-g(W+hW,X)hY{\big]}&\\
\nonumber&+\mu{\big[}g(W+hW,Y)hX-g(W+hW,X)hY&\\
\nonumber&+\{(1-k)g(W,X)-g(W,hX)+\eta(W)\eta(hX)\}\eta(Y)\xi&\\
\nonumber&-\{(1-k)g(W,X)-g(W,hY)+\eta(W)\eta(hY)\}\eta(X)\xi{\big]}&\\
\nonumber&-B(\phi W){\big[}\eta(Y)X-\eta(X)Y{\big]}-A(\phi W)[k\{\eta(Y)X-\eta(X)Y\}&\\
\nonumber&+\mu\{\eta(Y)h X-\eta(X)hY\}{\big]}.
\end{align}
This leads  to the following:
\begin{theorem}
  In a generalized $\phi$-recurrent generalized $(k,\mu)$-contact metric manifold, the curvature tensor $R$ satisfies the relation $(\ref{eqn4.17})$.
\end{theorem}
\section{example of generalized $\phi$-recurrent generalized $(k,\mu)$-contact metric manifold }
\begin{example}
We consider a $3$-dimensional manifold $M=\{x,y,z) \in \mathbb{R}^3:x \neq 0, y\neq 0\}$, where $(x,y,z)$ are the standard coordinates in $\mathbb{R}^3$.
Let $\{E_1,E_2,E_3\}$ be a linearly independent global frame on $M$ given by
  \begin{equation*}
    E_1=\frac{\partial}{\partial y},\ \ \ \ E_2=2xy\frac{\partial}{\partial z},\ \ \ \ \ E_3=\frac{\partial}{\partial z}.
  \end{equation*}
  Let $g$ be the Riemannian metric defined by $g(E_1,E_3)=g(E_2,E_3)=g(E_1,E_2)=0$, $g(E_1,E_1)=g(E_2,E_2)=g(E_3,E_3)=1.$
Let $\eta$ be the $1$-form defined by $\eta(U)=g(U,E_3)$ for any $U \in \chi(M)$. Let $\phi$ be the $(1,1)$ tensor field defined by $\phi E_1=-E_2$,
$\phi E_2=E_1$ and $\phi E_3 =0$. Then using the linearity of $\phi$ and $g$ we have $\eta(E_3)=1$, $\phi^2 U=-U+\eta(U)E_3$ and $g(\phi U,\phi W)=g(U,W)-
\phi(U)\phi(W)$ for any $U,W \in \chi(M)$. Moreover $hE_1=-E_1$, $hE_2=E_2$ and $h E_3=0$. Thus for $E_3=\xi$, $(\phi, \xi, \eta, g)$ defines a contact metric structure on $M$. Let $\nabla$ be the Riemannian connection of $g$. Then we have
$$[E_1,E_2]=\frac{1}{y}E_2,\ \ \ [E_1,E_3]=0,\ \ \ \ \ [E_2,E_3]=0.$$
Using Koszul  formula for the Riemannian metric $g$, we can easily calculate
\begin{align*}
  \nabla_{E_1} E_1=0, \ \ \ \ \nabla_{E_1}E_2=0, \ \ \ \ \ \nabla_ {E_1}E_3=0, &  \\
  \nabla_{E_2}E_1=-\frac{1}{y}E_2,\ \ \ \nabla_{E_2}E_2=\frac{1}{y}E_1, \ \ \ \ \nabla _{E_2}E_3=0,&  \\
   \nabla_{E_3}E_1=0,\ \ \ \nabla_{E_3}E_2=0, \ \ \ \ \nabla_{E_3}E_3=0. &
\end{align*}
From the above it can be easily seen that $(\phi,\xi,\eta,g)$ is a generalized $(k,\mu)$-contact metric structure on $M$.
Consequently $M^3(\phi,\xi,\eta,g)$ is a generalized $(k,\mu)$-contact metric manifold with $k=-\frac{1}{y}$ and $\mu=-\frac{1}{y}.$\\
\indent Using the above relations, we can easily calculate the non-vanishing components of the curvature tensor as follows:\\
$$R(E_1,E_2)E_1=\frac{2}{y^2}E_2,\ \ \   R(E_1,E_2)E_2=-\frac{2}{y^2}E_1,$$
and the components which can be obtained  from these by the symmetry properties. We shall now show that such a generalized $(k,\mu)$-contact metric manifold
is generalized $\phi$-recurrent. Since $ \{E_1,E_2,E_3\}$ forms a basis of $M^3$, any vector field $X,Y,Z\in \chi(M)$ can be written as
$$X=a_1E_1+b_1E_2+c_1E_3,$$
$$Y=a_2E_1+b_2E_2+c_2E_3,$$
$$Z=a_3E_1+b_3E_2+c_3E_3,$$
where $a_i,b_i,c_i\in {\mathbb{R}}^+$ (the set of all positive real numbers), $i=1,2,3.$ Then
\begin{equation}\label{eqn5.1}
  R(X,Y)Z=\frac{2}{y^2}(a_1 b_2-a_2 b_1)(a_3 E_2-b_3 E_1)
\end{equation}
and
\begin{eqnarray}
\label{eqn5.2}
  G(X,Y)Z &=& (a_2a_3+b_2b_3+c_2c_3)(a_1E_1+b_1E_2+c_1E_3)\\
  \nonumber &-&(a_1a_3+b_1b_3+c_1c_3)(a_2E_1+b_2E_2+c_2E_3).
  \end{eqnarray}
By virtue of (\ref{eqn5.1}) we have the following:
\begin{equation}\label{eqn5.3}
  (\nabla_{E_1}R)(X,Y)Z=\frac{4}{y^3}(a_1b_2-a_2b_1)(b_3E_1-a_3E_2),
\end{equation}
  \begin{equation}\label{eqn5.4}
    (\nabla_{E_2}R)(X,Y)Z=0,
  \end{equation}
  \begin{equation}\label{eqn5.5}
    (\nabla_{E_3}R){(X,Y)Z}= 0.
  \end{equation}
 From (\ref{eqn5.1}) and (\ref{eqn5.2}), we get
  \begin{equation}\label{eqn5.6}
    \phi^2(R(X,Y)Z)=u_1E_1+u_2E_2\ \ \ \ and \ \ \ \ \ \phi^2(G(X,Y)Z)=v_1E_1+v_2E_2,
  \end{equation}
 where
 \begin{equation*}
   u_1=\frac{2b_3}{y^2}(a_1b_2-a_2b_1),\ \ \ \ \ u_2=-\frac{2a_3}{y^2}(a_1b_2-a_2b_1),
 \end{equation*}
 \begin{equation*}
   v_1=a_2(b_1b_3+c_1c_3)-a_1(b_2b_3+c_2c_3),
 \end{equation*}
  \begin{equation*}
    v_2=b_2(a_1a_3+c_1c_3)-b_1(a_2a_3+c_2c_3).
  \end{equation*}
  Also from (\ref{eqn5.3})-(\ref{eqn5.5}), we obtain
  \begin{equation}\label{eqn5.7}
    \phi^2((\nabla_{E_i}R)(X,Y)Z)=p_iE_1+q_iE_2 \ \ \ \ \ i=1,2,3,
  \end{equation}
  where
  \begin{equation*}
 p_1=-\frac{4b_3}{y^3}(a_1b_2-a_2b_1),\ \ \ q_1=\frac{4a_3}{y^3}(a_1b_2-a_2b_1),
  \end{equation*}
  \begin{equation*}
 p_2=0,\ \ q_2=0,\ \ p_3=0,\ \ \ \ q_3=0.
  \end{equation*}
  Let us consider the $1$-forms as
  \begin{align}\label{eqn5.8}
 A(E_1)=\frac{v_2p_1-v_1q_1}{u_1v_2-u_2v_1},\ \ \ \ \ & B(E_1)=\frac{u_1q_1-u_2p_1}{u_1v_2-u_2v_1}, \\
\nonumber A(E_2)=0,\ \ \ \ & B(E_2)=0, \\
\nonumber A(E_3)=0,\ \ \ \ & B(E_3)=0,
  \end{align}
  where $a_1b_2-a_2b_1 \neq 0,$ $v_2p_1-v_1q_1 \neq 0$, $u_1q_1-u_2p_1 \neq 0$, $u_1v_2-u_2v_1 \neq 0$.
\end{example}
From (\ref{eqn4.1}), we have
\begin{equation}\label{eqn5.9}
  \phi^2((\nabla_{E_i} R)(X,Y)Z)=A(E_i)\phi^2(R(X,Y)Z)+B(E_i)\phi^2(G(X,Y)Z),\ \ \ i=1,2,3.
\end{equation}
By virtue of (\ref{eqn5.6})-(\ref{eqn5.8}), it can be easily shown that the manifold satisfies the relation (\ref{eqn5.9}).
Hence the manifold under consideration is a $3$-dimensional generalized $\phi$-recurrent generalized $(k,\mu)$-contact metric manifold which is neither
$\phi$-symmetric nor $\phi$-recurrent.\\
This leads the following:
\begin{theorem}
  There exists a $3$-dimensional generalized $\phi$-recurrent generalized $(k,\mu)$-contact metric manifold, which is neither $\phi$-symmetric nor
  $\phi$-recurrent.
\end{theorem}

\vspace{0.1in}
\noindent Shyamal Kumar Hui\\
Department of Mathematics\\
The University of Burdwan \\
Burdwan, 713104\\
 West Bengal, India\\
E-mail: skhui@math.buruniv.ac.in


\begin{thebibliography}{15}
\bibitem{ARSL}
Arslan, K., De, U. C., Murathan, C. and Yildiz, A., \emph{On generalized recurernt Riemannian manifollds}, Acta Math. Hungarica, {\bf 123} (2009), 27-39.
\bibitem{BAISHYA}
Baishya, K. K., Eyasmin, S. and Shaikh, A. A., \emph{On $\phi$-recurrent generalized $(k, \mu)$-contact metric manifolds}, Lobachevskii J. of Math., {\bf 27} (2007), 3-13.
\bibitem{BLAIR1}
Blair, D. E., \emph{Contact manifolds in Riemannian geometry}, Lecture Notes in Math. {\bf 509}, Springer-Verlag, 1976.
\bibitem{BLAIR2}
Blair, D. E., \emph{Two remarks on contact metric structures}, Tohoku Math. J., {\bf 29} (1977), 319-324.
\bibitem{BLAIR3}
Blair, D. E., Koufogiorgos, T. and Papantoniou, B. J., \emph{Contact metric manifolds satisfying a nullity condition}, Israel J. of Math., {\bf 19}
(1995), 189-214.
\bibitem{BOEC}
Boeckx, E., \emph{A full classification of contact metric $(k,\mu)$-spaces}, Illinois J. Math., {\bf 44} (2000), 212-219.
\bibitem{KOUF}
Koufgiorogs, T. and Tsichlias, C., \emph{On the existence of new class of contact metric manifolds}, Canad. Math. Bull., {\bf XX(Y)} (2000), 1-8.
\bibitem{DE}
De, U. C. and Guha, N., \emph{On generalized recurrent manifolds}, J. Nat. Acad. Math. India, {\bf 9} (1991), 85-92.
\bibitem{DE1}
De, U. C. and Guha, N. and Kamilya, D., \emph{On generalized Ricci recurrent manifolds}, Tensor, N. S., {\bf 56} (1995), 312-317.
\bibitem{DUBEY}
Dubey, R. S. D., \emph{Genralized recurrent spaces}, Indian J. Pure, Appl. Math., {\bf 10} (1979), 1508-1513.
\bibitem{PATT}
Patterson, E. M., \emph{Some theorems on Ricci recurrent spaces}, J. London Math. Soc., {\bf 27} (1952), 287-295.
\bibitem{SHAIKH3}
Shaikh, A. A. and Ahmad, H., \emph{On generalized $\phi$-recurrent Sasakian manifolds}, Applied Mathematics, {\bf 2} (2011), 1317-1322.
\bibitem{SHAIKH1}
Shaikh, A. A. and Baishya, K. K., \emph{On a contact metric manifold}, Diff. Geom. Dynamical System, {\bf 8} (2006), 253-261.
\bibitem{SHAIKH2}
Shaikh, A. A., Arslan, K., Murathan, C. and Baishya, K. K., \emph{On $3$-dimensional generalized $(k,\mu)$-contact metric manifold}, Balkan J. of Geom. and its Applications, {\bf 12} (2007), 122-134.
\bibitem{TAKA}
Takashashi, T., \emph{Sasakian $\phi$-symmetric spaces}, Tohoku Math. J., {\bf 29} (1977), 91-113.
\bibitem{TANNO}
Tanno, S., \emph{Ricci curvatures of contact Riemannian manifolds}, Tohoku Math. J., {\bf 40} (1988), 441-448.
\bibitem{WALK}
Walker, A. G., \emph{On Ruse's spaces of recurrent curvature}, Proc. London Math. Soc., {\bf 52} (1950), 36-64.
\end{thebibliography}
\end{document}